\documentclass{amsart}
\usepackage{amssymb,amsmath,amscd,enumerate}
\usepackage[all]{xy}
\usepackage{hyperref}

\newtheorem{thm}{Theorem}[section]

\newtheorem{corollary}[thm]{Corollary}
\newtheorem{prop}[thm]{Proposition}
\newtheorem{lemma}[thm]{Lemma}
\newtheorem{fact}[thm]{Fact}

\theoremstyle{definition}
\newtheorem{defn}[thm]{Definition}
\newtheorem{example}[thm]{Example}

\theoremstyle{remark}
\newtheorem{remark}[thm]{Remark}

\newcommand{\bt}{\begin{thm}}
\newcommand{\et}{\end{thm}}
\newcommand{\bp}{\begin{prop}}
\newcommand{\ep}{\end{prop}}
\newcommand{\bd}{\begin{defn}}
\newcommand{\ed}{\end{defn}}
\newcommand{\bl}{\begin{lemma}}
\newcommand{\el}{\end{lemma}}
\newcommand{\bfa}{\begin{fact}}
\newcommand{\efa}{\end{fact}}
\newcommand{\bc}{\begin{corollary}}
\newcommand{\ec}{\end{corollary}}
\newcommand{\bex}{\begin{example}}
\newcommand{\eex}{\end{example}}
\newcommand{\br}{\begin{remark}}
\newcommand{\er}{\end{remark}}
\newcommand{\ben}{\begin{enumerate}}
\newcommand{\een}{\end{enumerate}}

\newcommand{\mf}{\mathbf}
\newcommand{\llength}{\mbox{\em length}}
\newcommand{\length}{\mbox{length}}
\newcommand{\mx}{\mathfrak{m}}

\newcommand{\sotto}[2]{#1_{#2}}

\newcommand{\rrr}{\rightarrow}
\newcommand{\ra}{\rightarrow}

\newcommand{\ds}{\displaystyle}

\newcommand{\ideal}[1]{\sotto {{\mathcal I}}{#1}}

\newcommand{\exact}[3]
{0 \rrr #1 \rrr #2
\rrr #3 \rrr 0}

\newcommand{\pso}{\mathbb{P}^3}

\newcommand{\PP}{\mathbb{P}}

\newcommand{\Z}{\mathbb{Z}}

\newcommand{\N}{\mathbb{N}}

\newcommand{\coo}{{\mathcal O}}

\begin{document}

\title{Smooth curves specialize to extremal curves}

\author[R. Hartshorne]{Robin Hartshorne}
\address{Department of Mathematics, University  of
California, Berkeley, California 94720--3840}
\email{robin@math.berkeley.edu}

\author[P.~Lella]{Paolo Lella}
\address{Dipartimento di Matematica dell'Universit\`{a} di Torino\\
         Via Carlo Alberto 10,
         10123 Torino, Italy}
\email{\href{mailto:paolo.lella@unito.it}{paolo.lella@unito.it}}
\urladdr{\href{http://www.dm.unito.it/dottorato/dottorandi/lella/}{http://www.dm.unito.it/dottorato/dottorandi/lella/}}

\author[E.~Schlesinger]{Enrico Schlesinger}
\address{Dipartimento di Matematica del Politecnico di Milano\\
         Piazza Leonardo da Vinci 32,
         20133 Milano, Italy}
\email{\href{mailto:enrico.schlesinger@polimi.it}{enrico.schlesinger@polimi.it}}
\urladdr{\href{http://www1.mate.polimi.it/~enrsch//}{http://www1.mate.polimi.it/~enrsch/}}

\thanks{
The third author was partially supported by
MIUR PRIN 2007:
  {\em 	 Moduli, strutture geometriche e loro applicazioni}.
 }

\subjclass[2000]{14C05, 14H50, 14Q05, 13P10}
\keywords{Hilbert scheme, locally Cohen-Macaulay curve, initial ideal, weight vector, Groebner bases, smooth curves}

\begin{abstract}
Let $H_{d,g}$ denote the Hilbert scheme of locally Cohen-Macaulay curves of degree $d$ and genus $g$
in projective three space. We show that, given a smooth irreducible curve $C$ of degree $d$ and genus $g$,
there is a rational curve $\{[C_t]: t \in \mathbb{A}^1\}$ in $H_{d,g}$ such that $C_t$ for $t \neq 0$
is projectively equivalent to $C$, while the special fibre $C_0$ is an extremal curve. It follows
that smooth curves lie in a unique connected component of $H_{d,g}$. We also determine necessary and sufficient conditions
for a locally Cohen-Macaulay curve to admit such a specialization to an extremal curve.
\end{abstract}

\maketitle
%
%
%
%
%
%

\section{Introduction}
In this paper we study curves in projective three space $\PP^3$ over an algebraically closed field $k$.
While smooth curves are often the first to be studied, the development of liaison theory in recent years
\cite{MDP} has shown that these are best understood in the context of the larger class of locally Cohen-Macaulay
curves, that is, one-dimensional closed subschemes of $\PP^3$ with no isolated points and no embedded points.
For each $d$ and (arithmetic) genus $g$, the locally Cohen-Macaulay curves form an open subset $H_{d,g}$ of the
full Hilbert scheme $\mbox{Hilb}_{d,g}$ of all closed subschemes of $\PP^3$ with Hilbert polynomial $dn+1-g$.
The smooth curves form an open subset $H^0_{d,g}$ of $H_{d,g}$.

It is known from classical examples that $H^0_{d,g}$ may not be connected. On the other hand, the thesis of the first author
\cite{hthesis} showed that for any $d,g$ the full Hilbert scheme $\mbox{Hilb}_{d,g}$ is connected whenever it is nonempty.
So it is natural to ask whether the Hilbert scheme $H_{d,g}$ of locally Cohen-Macaulay curves is connected, a question first raised
by Martin-Deschamps and Perrin in their paper \cite{extremal}. Up to now, the connectedness of $H_{d,g}$ has been established
only for very small degree $d$, or for very large genus $g$; more precisely, for $d \leq 4$   \cite{nthree,ns-comp} and
$g \geq \frac{1}{2} (d\!-\!3)(d\!-\!4)-1$  \cite{AA1,hmdp3,sabadini}). The method is to identify all the irreducible components of $H_{d,g}$
and then to verify that they intersect each other in such a way as to make the whole set connected. It is the problem of
finding all irreducible components of $H_{d,g}$ that seems to block any further progress in this direction.

Another approach makes use of the so-called {\em extremal curves}. Martin-Deschamps and Perrin showed in two papers \cite{bounds,extremal}
that for any curve $C$ of degree $d$ and genus $g$ the dimensions of the cohomology modules $H^{1} (\PP^3,\ideal{C}(n))$
are bounded by explicit functions of $d,g$ (see exact statement below). Those curves that attain the maximum are called
{\em extremal curves}, and the family of extremal curves forms a non-empty closed irreducible component of $H_{d,g}$
for each $d,g$. Thus one can ask which curves can be connected, through a chain of irreducible components of $H_{d,g}$
to an extremal curve. Various methods have shown for example that this is possible for smooth rational or elliptic curves,
for ACM curves, for any curve whose Rao module is a complete intersection, or for any curve in the biliaison class
of an extremal curve \cite{hconn,footnote,perrin-pas}.

The purpose of this paper is to give a necessary and sufficient condition for an irreducible component $V$ of $H_{d,g}$
to contain an extremal curve. We say that a curve $C$ in $\PP^3$ satisfies condition (*) if there is a point $p$ in $\PP^3$
such that:

\noindent
\begin{itemize}
  \item[1)]  for every line $N$ containing $p$
  the degree of $C \cap N$ is at most two; furthermore, $\deg(C \cap N)= 2$  for at most finitely
  many such lines;
  \item[2)]  if $N$ is a line through $p$ meeting $C$ in a scheme of degree $2$,
then there is a plane $K$ through $N$ such that the divisor $2N$ in $K$
meets $C$ in a scheme of degree $\leq 3$.
\end{itemize}
Note that smooth irreducible curves satisfy condition (*).
Our main result shows that an irreducible component $V$ of $H_{d,g}$ contains an extremal curve if and only if
a general curve $C \in V$ satisfies (*). We prove also the stronger result that if $C$ satisfies (*), then $C$ {\em specializes}
to an extremal curve, meaning that there is a flat family $\{C_t\}_{t \in \mathbb{A}^1}$ such that $C=C_1$, for all $t \neq 0$,
$C_t$ is obtained from $C$ by an automorphism of $\PP^3$, and $C_0$ is an extremal curve. In the last section of the paper
we give a geometric characterization of curves satisfying condition (*). From this it follows that in characteristic zero
a reduced curve that has embedding dimension $\leq 2$ at every point satisfies condition (*).

As a corollary of our theorem, we recover many earlier results, and in particular we show that any smooth irreducible curve specializes to an extremal
curve. So for example we can conclude that the two irreducible components of $H^0_{9,10}$, neither of which meets the closure of the other
 \cite[IV,6.4.3]{AG}, can be connected via extremal curves in $H_{9, 10}$.

The method we use, which appears already in the earlier paper \cite{ls}  of the second and third authors, is to place the
curve on a certain surface of the kind called a monoid by Cayley in his early study of space curves \cite{cayley}. The limit of the curve under a weighted automorphism
of $\PP^3$ is then a multiple line on the limit surface, which one shows to be an extremal curve. The essential point, of course, is
to show that the limit has no embedded points, and this is accomplished by a careful control of the degree and genus
of the curves involved.

 We should mention also that there are curves to which our results do not apply. For $d=4$ and genus $\leq -3$,
there is a closed irreducible component of $H_{d,g}$, made up of the so-called {\em thick} $4$-lines,
that does not meet the component of extremal curves \cite{ns-comp}.

\vspace{0.5cm} \noindent {\bf Acknowledgments.}
We thank the American Institute of Mathematics  and the organizers of the AIM workshop {\em Components of Hilbert Schemes} Robin Hartshorne, Diane Maclagan, and Gregory Smith \cite{aim}.  The seeds of this paper were planted there. The  discovery of the results of this paper was
accomplished with the help of several computations performed with {\em Macaulay2} \cite{M2} using the package {\em gfanInterface}.
We would like to thank D. Grayson and M. Stillman for {\em Macaulay 2} and A. Jenssen for the program Gfan \cite{gfan}.
The first author would like to add that the method of specialization used here is entirely the discovery of the second and third authors;
his contribution was to help find the best result that could be obtained by that method. We would like to thank the referee for pointing out an error in the original proof of Proposition 2.11, and for suggesting an argument to correct it.

\section{Construction of specializations to extremal curves} \label{two}
In this section we establish notation and terminology and review some known results that we will need later.
We work over an algebraically closed field $k$ of arbitrary characteristic. We denote by the symbol
$\ideal{X}$ the ideal sheaf of a subscheme $X \subset \PP^3$. Given a coherent
sheaf $\mathcal{F}$ on $\PP^3$, we define $h^i(\mathcal{F})= \dim \, H^{i}(\PP^3,\mathcal{F})$ and
$\ds H^i_* (\mathcal{F})= \bigoplus_{n \in \Z} H^{i}(\PP^3,\mathcal{F}(n))$. We write
$R=k[x,y,z,w]$ for the homogeneous coordinate ring $H^{0}_* (\coo_{\PP^3})$ of $\PP^3$.

\bd
A {\em curve} in $\PP^3$ (or more precisely a locally Cohen-Macaulay curve)
is a one dimensional subscheme $C \subset \PP^3$ without zero dimensional associated points; this means
that all irreducible components of $C$ have dimension $1$, and that $C$ has no embedded points.
\ed
We denote by $H_{d,g}$ the Hilbert scheme parametrizing curves of degree $d$ and arithmetic
genus $g$ in $\PP^3$. This is an open subscheme of the full Hilbert scheme parametrizing
all one dimensional subschemes of $\PP^3$ with Hilbert polynomial $dn+1-g$. Martin-Deschamps and Perrin \cite{bounds}
have found sharp bounds for the Rao function $h^1 (\ideal{C}(n))$ of a curve in $H_{d,g}$. To state these bounds we need
to introduce some notation. Let $a=\frac{1}{2}  (d\!-\!2)(d\!-\!3)-g$ and $l=d-2$. Then for
$d \geq 2$ and $g \leq \frac{1}{2}  (d\!-\!2)(d\!-\!3)$ define
$$
\rho_{d,g} (n) =
\left\{
\begin{array}{ll}
0 &
\mbox{ if \;\;\; $  n \leq -a$}
\\
n+a &
\mbox{ if \;\;\; $ -a \leq n \leq 0$}
\\
a
&
\mbox{ if \;\;\; $ 0 \leq n \leq l$}
\\
a+l-n
&
\mbox{ if \;\;\; $ l \leq n \leq a+l$}
\\
0 &
\mbox{ if \;\;\; $  n \geq a+l$}
\end{array}
\right.
$$
\bt[\cite{bounds,subextremal}] \label{bounds}
Let $C \subset \pso$ be a curve of degree $d$ and arithmetic
genus $g$. Then
\begin{enumerate}
\item
$C$ is a plane curve if and only if $g = \frac{1}{2}  (d\!-\!1)(d\!-\!2)$; in this case
$h^1 (\ideal{C}(n))=0$ for every $n \in \Z$.
\item
If $C$ is not a plane curve, then $d \geq 2$, $g \leq \frac{1}{2}  (d\!-\!2)(d\!-\!3)$
and
$$
h^1 (\ideal{C}(n)) \leq \rho_{d,g} (n) \;\;\; \mbox{for all $n \in \Z$.}
$$
\end{enumerate}
\et

\br
The set of pairs for which $H_{d,g}$ is non-empty was determined by Hartshorne in \cite{hgenus}, see also
\cite{ok-sp}.
The bound in the theorem is proven in \cite{bounds} in characteristic zero; Nollet in \cite{subextremal} has shown
that the same result is valid in characteristic $p$ as well.
\er

\bd
A curve $E \subset \pso$ of degree $d$ and genus $g$ is called {\em
extremal}  if either it is a plane curve or
it satisfies  $h^1 (\ideal{E}(n)) = \rho_{d,g}(n)$ for every $n \in \Z$.
\ed
\br
Our definition of extremal curves is equivalent to the one given by Hartshorne \cite{hconn}.
Martin-Deschamps and Perrin required extremal curves not to be ACM, that is,
$g<\frac{1}{2} (d\!-\!2)(d\!-\!3)$.
\er
It is a remarkable discovery of Martin-Deschamps and Perrin  that extremal curves exist for every $d,g$ for which $H_{d,g}$
is non-empty:

\bt[\cite{hgenus,bounds,extremal}]\label{mdp-extremal}
If $H_{d,g}$ is non-empty, then there exist
extremal curves of degree $d$ and genus $g$, and they form a closed irreducible component $E_{d,g}$
of $H_{d,g}$.
\et

\br
In fact in \cite{bounds} it is shown that, for  given $d$ and $g$,  extremal curves in $H_{d,g}$ have
the same  Hilbert function, that is, the functions $h^0 (\ideal{E}(n))$ and $h^2 (\ideal{E}(n))$
are constant on $E_{d,g}$; furthermore,
\begin{equation}\label{extr_def}
h^{i} (\ideal{C}(n)) \leqslant h^{i} (\PP^3,\ideal{E}(n)) \qquad \mbox{for every $ n \in \Z$ and for $i=0,1,2$}
\end{equation}
for every curve $C$ in $H_{d,g}$ and every extremal curve $E \in E_{d,g}$.
Thus extremal curves have the largest possible cohomology.
\er

\br
Suppose $d \geq 2$ and $g < \frac{1}{2}  (d\!-\!2)(d\!-\!3)$, that is, $a>0$ and $l \geq 0$.
Martin-Deschamps and Perrin show that
the extremal curves are precisely the minimal curves associated to {\em extremal Koszul modules}. These are
the modules of the form
$ R/(l_1,l_2,F,G)$
where $l_1$ and $l_2$ are linear forms, $F$ and $G$ are forms of degree $a$ and $a+l$
respectively, and $l_1,l_2,F,G$ have no common zeros. A curve is extremal if and only if
its Rao module is an extremal module shifted to the left by $a\!-\!1$.
\er

The following proposition describes extremal curves supported on the line $x=y=0$.

\bp[\cite{ls} Proposition 2.5] \label{extr_mls}
Let $(d,g)$ be a pair of integers satisfying $d \geq 2$ and $g < \frac{1}{2} (d\!-\!2)(d\!-\!3)$. Let $F$ and $G$
be two forms in $k[z,w]$ of degrees $a$ and $a+l$, respectively, with no common zeros.
The surface $S$ of equation $xG-y^{d-1}F=0$ is irreducible and generically smooth along the line
$L$ of equations $x=y=0$. It therefore contains a unique curve $E$ of degree $d$ supported on $L$.
The curve $E$ is extremal of degree $d$ and genus $g$, and its Rao module is
$$
H^{1}_* (\ideal{E}) \cong R/(x,y,F,G) (a\!-\!1) \cong k[z,w]/(F,G) (a\!-\!1)
$$
The homogeneous ideal of $E$ is generated by $x^2$, $xy$, $y^d$ and $xG-y^{d-1}F$.
\ep

\begin{proof}
The surface $S$ is irreducible because $F$ and $G$ have no common zeros, and it is smooth at points
of $L$ where $G$ is different from zero. Therefore  the ideal of $L$ in the local ring
$\coo_{S,\xi}$ of the generic
point $\xi$ of $L$ is generated by one function $t$ , and the ideal of a curve of degree $d$ supported on $L$
must be $t^d\coo_{S,\xi}$. Since a locally Cohen-Macaulay curve supported on $L$ is determined by its ideal
at the generic point of $L$, we see that there is a unique curve $D_m \subset S$ supported on $L$ of degree $m$
for every $m \geqslant 1$. For $m=d-1$, the curve $P=D_{d-1}$ is the planar multiple structure of equations
$x=y^{d-1}=0$. 
We note that
$\ideal{P} \otimes \coo_L \cong \coo_L (-1) \oplus \coo_L (1-d)$  where the two generators are the images
of $x$ and $y^{d-1}$. The two forms $F$ and $G$ define a surjective map $\coo_L (-1) \oplus \coo_L (1-d) \ra \coo_L (a\!-\!1)$;
composing this with the natural map  $\ideal{P} \ra \ideal{P} \otimes \coo_L$  we obtain
a surjection $\phi:\ideal{P} \ra \coo_L (a\!-\!1)$. We let $E$ be the subscheme of $\PP^3$ whose ideal sheaf is
the kernel of $\phi$. By construction we have an exact sequence
 \begin{equation}\label{extr_seq}
    \exact{\ideal{E}}{\ideal{P}}{\coo_L (a\!-\!1)}
 \end{equation}
This sequence shows that $E$ is a (locally Cohen-Macaulay) curve of degree $d$ and genus $g$, and that its homogeneous
ideal is generated by $x^2$,$xy$,$y^d$ and $xG-y^{d-1}F$. Therefore $E=D_d$ is the unique curve of degree $d$ contained
in $S$ and supported on $L$. Finally, the long exact cohomology sequence associated to (\ref{extr_seq}) shows that the Rao
module of $E$ is
$$
M_E= k[z,w] /(F,G) (a\!-\!1) = R/(x,y,F,G) (a\!-\!1).
$$
Hence $E$ is an extremal curve.
\end{proof}

\bigskip

Given a weight vector $\omega \in \N^{N+1}$, the $\omega$-degree of
a monomial $\boldsymbol{x}^{\mathbf{n}} = x_0^{n_0}\cdots x_N^{n_N}$ is defined as $\deg_\omega \boldsymbol{x}^{\mathbf{n}} = \mathbf{n}\cdot\omega = n_0\omega_0 + \ldots + n_N\omega_N$, and the $\omega$-degree of an arbitrary polynomial
$P(\boldsymbol{x})=\sum_{\mf{n}}c_{\mathbf{n}} \boldsymbol{x}^\mathbf{n}$ is the maximum degree of a monomial appearing in $P(\boldsymbol{x})$. In analogy with the notion of leading term of a polynomial with respect to a term ordering, one defines the \emph{initial form} of $P(\boldsymbol{x})$ with respect to $\omega$ as
\[
\text{in}_\omega\big(P(\boldsymbol{x})\big) = \sum_{\begin{subarray}{c} \mathbf{n}\\ \deg_{\omega} \boldsymbol{x}^{\mathbf{n}} = \deg_{\omega} P\end{subarray}} c_{\mathbf{n}}\boldsymbol{x}^{\mathbf{n}}.
\]
For any ideal $I \subset k[x_0,\ldots,x_N]$ one defines the \emph{initial ideal} $\text{in}_\omega (I) $ of $I$ with respect to  $\omega$ as the ideal generated
by all initial forms $ \text{in}_{\omega} (P)$ as $P$ varies in $I$. We will use the following well known fact
(see for example \cite{bayer-m}, \cite[Theorem 15.17]{eisenbud} or \cite[Theorem 4.3.22]{KreuzerRobbiano2}).

\bp \label{prop:flatWeight}
Given a subscheme $X \subset \PP^N$ and a weight vector $\omega \in \N^{N+1}$, there is a flat family $X_t \subset \PP^N \times \mathbb{A}^1$,
whose fibres $X_t$ for $t \neq 0$
are isomorphic to $X=X_1$ via an automorphism of $\PP^N$ and whose special fibre $X_0$ is the subscheme of $\PP^N$ defined
by the initial ideal of $X$ with respect to $\omega$.
\ep

%

\bigskip

Our basic  remark is that the homogeneous polynomial $xG-y^{d-1}F$, if $F$ and $G$ are forms in $z$ and $w$,
is also homogeneous with respect to the weight vector $\omega=(d,2,1,1)$. It follows that the
ideal of the extremal curve $E$ of Proposition \ref{extr_mls} is homogeneous with respect to the grading
defined by $\omega$, and so it coincides with its initial ideal with respect to $\omega$;
we can then construct flat families having $E$ as a special fibre by taking the initial
ideal with respect to this particular weight vector.

\bp \label{method}
Let $(d,g)$ be a pair of integers satisfying $d \geq 2$ and $g < \frac{1}{2} (d\!-\!2)(d\!-\!3)$.
Let $\omega=(d,2,1,1)$.
Let $S$ be a surface whose defining  equation $f$ satisfies
\[
\textnormal{in}_\omega (f)=xG-y^{d-1}F
\]
where $F$ and $G$
be two forms in $k[z,w]$ of degrees $a$ and $a+l$, respectively, with no common zeros.
If $C$ is a curve of degree $d$ and genus $g$ on $S$ that does not meet the line $z=w=0$, then
the saturation of the initial ideal of $C$ with respect to $\omega$ is in the ideal
generated by $x^2$, $xy$, $y^d$ and $xG-y^{d-1}F$ and therefore defines an extremal curve.
\ep

\begin{proof}
Suppose $C$ is a curve that does not meet the line $M$ of equation $z=w=0$.
The projection of $C$ from the point $P=[1:0:0:0]$ on the plane $x=0$ is a curve
$\Gamma$ of some degree $m \leq d$, whose equation in the plane is of the form
$$
g=g_0 y^m +g_1(z,w) y^{m-1}+ \cdots + g_m(z,w)=0
$$
where $g_{i}(z,w)$ is a form of degree $i$ in $z$ and $w$. Since $C$ does not meet the line $z=w=0$,
it follows that the point $[0:1:0:0]$ is not on the curve $\Gamma$, so $g_0$ is non-zero, and
$\text{in}_\omega (g)= g_0 y^m$. Now $C$ lies
on the cone over $\Gamma$, which has the same equation $g=0$.  Therefore
$y^m$ belongs to $\text{in}_\omega (I_C)$.
Similarly, projecting from $[0:1:0:0]$ to the plane $y = 0$ shows that some power of $x$ belongs to $\text{in}_\omega (I_C)$.
We conclude that the curve $C_0$ defined by $\text{in}_\omega (I_C)$ is supported on the line $L$ of equation $x=y=0$.

Now we add the hypothesis that $C$ is contained in a surface $S$ with an equation $f$ such that
$\text{in}_\omega (f)=xG-y^{d-1}F$. Then the scheme
$C_0$ is contained in the surface $S_0$ of equation $xG-y^{d-1}F=0$. By flatness, the Hilbert polynomial of $C_0$ coincides
with that of $C$, so $C_0$ is a one dimensional subscheme of $\PP^3$ of degree $d$ and genus $g$.
Let $E$ be the largest Cohen-Macaulay curve contained in $C_0$: it is the curve of degree $d$ obtained
from $C_0$ throwing away its embedded points. By Proposition \ref{extr_mls} $E$ is the unique
curve of degree $d$ contained in $S_0$ and supported on the line $L$; it is the extremal curve
whose ideal
is generated by $x^2$, $xy$, $y^d$ and $xG-y^{d-1}F$. Since $E \subset C_0$ and the two schemes have the same Hilbert polynomial, we conclude
$E=C_0$, and the statement is proven.
\end{proof}

\section{Specialization of curves satisfying condition (*)} \label{three}

Now we come to our main result.
We denote by  $X \cap Y$ the scheme theoretic intersection of two closed subschemes
of an ambient scheme $V$. If $Z$ is a zero dimensional scheme of finite type over $k$
and $p$ is a closed point in the support of $Z$,
we denote by $\textnormal{mult}_p (Z)=\length (\coo_{Z,p})$ the multiplicity of $Z$ at $p$, and by
$\ds \deg(Z)= \sum_{p \in \textnormal{Supp} (Z)} \textnormal{mult}_p (Z)$ the degree of $Z$.
For a curve $C$ in $\PP^3$ we consider a condition (*) that expresses the idea that $C$ behaves in some
ways like a smooth curve:

\vspace{.5cm}
\noindent
{\bf Condition (*)}

We say that a curve $C$ in $\PP^3$ satisfies condition (*) if there is a point $p$ in $\PP^3$
such that:

\noindent
\begin{itemize}
  \item[1)]  for every line $N$ containing $p$
  the degree of $C \cap N$ is at most two; furthermore, $\deg(C \cap N)= 2$  for at most finitely
  many such lines;
  \item[2)]  if $N$ is a line through $p$ meeting $C$ in a scheme of degree $2$,
then there is a plane $K$ through $N$ such that the divisor $2N$ in $K$
meets $C$ in a scheme of degree $\leq 3$.
\end{itemize}

\br
The set of points $p$ with respect to which $C$ satisfies condition (*) is open in $\PP^3$. If $2)$ holds, then, for a general
plane $K$ in the pencil of planes through $N$, the divisor $2N$ in $K$
meets $C$ in a scheme of degree $\leq 3$, and therefore, for at most finitely many planes
$K$ in the pencil, the divisor $2N$ in $K$  meets $C$ in a scheme of degree $\geq 4$.
\er

%

\bp \label{open-condition}
Condition (*) is an open condition on the Hilbert scheme.
\ep
\begin{proof}
Suppose that $C_0$ is a curve satisfying (*). It will be sufficient to show that for any flat family $\mathcal{C} = \{C_t\}_{t \in T}$ of curves, parametrized by an irreducible curve $T$, with $C_{t_0} = C_0$, a general curve $C_t$ in the family also satisfies (*). Choose $p \in \PP^3$ that satisfies condition (*) for $C_0$.
If a general curve in the family $\mathcal{C}$ does not satisfy (*), then either $a)$ for a
general $C_t$ all lines from $p$ to a point of $C_t$ meet $C_t$ at least twice, or $b)$ for a general $C_t$ there is a line $N$ through $p$ meeting $C_t$ at least $3$ times, or $c)$ there is a line $N$ through $p$ such that all planar double lines $2N$ containing $N$ intersect $C_t$ with multiplicity at least $4$. In each case, then, one can choose families of lines $N_t$ for each $C_t$ (possibly after a base extension $T'\rightarrow T$ of the family) whose limits
$N_0$ by semicontinuity will contradict the property (*) that holds for $C_0$. Hence a general $C_t$ in $\mathcal{C}$ must satisfy (*).
\end{proof}

\bt \label{main-result}
\
\begin{enumerate}[a)]
\item\label{it:mainThm_a} An irreducible component $V$ of $H_{d,g}$ contains an extremal curve if and only if
a general curve $C$ in $V$ satisfies condition (*).
\item\label{it:mainThm_b} If $C$ is a curve satisfying (*), then there is a flat family $\mathcal{C}=\{C_t\}_{t \in \mathbb{A}^1}$
such that $C_1=C$, for $t \neq 0$, $C_t$ is obtained from $C$ by an automorphism of $\PP^3$, and $C_0$
is an extremal curve (In this case we say $C$ {\em specializes} to an extremal curve).
\end{enumerate}
\et
\begin{proof}[Proof of Theorem \ref{main-result} (first part)]
\-

\noindent \textbf{Step 1.} First we deal with two special cases. If $g=\frac{1}{2} (d\!-\!1)(d\!-\!2)$,
then $H_{d,g}$ is irreducible and consists of the plane curves of degree $d$.
These obviously satisfy (*).
If $g=\frac{1}{2} (d\!-\!2)(d\!-\!3)$, again
$H_{d,g}$ is irreducible, and consists of ACM extremal curves. These have been described in \cite[3.5]{hgenus}.
The general such curve is a nonsingular twisted cubic for $d=3$, a nonsingular elliptic quartic for $d=4$, and for
$d \geq 5$ consists of a plane curve of degree $d\!-\!1$ with a line attached at one point. These general curves
clearly satisfy (*), so part a) of the theorem is true in these two cases.
For part $b)$ we can take the trivial family $C_t=C$ for all $t$. Thus for the remainder of this proof we may assume $g <\frac{1}{2} (d\!-\!2)(d\!-\!3)$ since these are the only remaining values where $H_{d,g}$ is non-empty.
(For $g=\frac{1}{2} (d\!-\!2)(d\!-\!3)$ it is not
true that every extremal curve satisfies (*): consider for example the curve defined by the ideal $(x^2,xy,y^{d-1})$ for
$d \geq 3$).

\noindent\textbf{Step 2.}
We will show next  that every extremal curve $C$ with $g <\frac{1}{2} (d\!-\!2)(d\!-\!3)$ satisfies (*).
We will use the classification of extremal curves of  \cite[Proposition 0.6]{extremal}. If the
curve $C$ 
is the disjoint union of a line and a plane curve, then $C$ obviously satisfies
(*). The ideal of any other extremal curve  can be brought with a change of variables into
the form
\[
I_C=(x^2,xy,Py^{b+2},xG+y^{b+1}FP)
\]
where $F,G$ are forms in $z,w$ of degrees $a$, $a+l$ respectively, with no common zeros, while
$P$ is a form in $y,z,w$ of degree $l-b$ that is not a multiple of $y$.
Here
$a= \frac{1}{2} (d\!-\!2)(d\!-\!3)-g \geq 1$, $l=d-2 \geq 1$, and $b \geq 0$.
The curve $C$ is contained in the double plane $2H=V(x^2)$, and, except
at points of the line $L=V(x,y)$, the curve is contained in the reduced plane $H=V(x)$
\cite[loc. cit.]{extremal}.

We claim that $C$ satisfies (*) with respect to the point $p=(1,0,0,0)$.
Since $C$ is contained in the double plane $2H$,
the intersection of $C$ with any line $N$ through $p$ has at most of degree $2$ in any case.
If $N$ is a line joining $p$ with a point $q$ of $C \setminus L$, then $C \cap N$ has degree one
because it is contained in  $C \cap H$. Thus we need to worry only about lines $N=pq$ where
$q=(0,0,z_0,w_0)$ is a point of $L$. The ideal of the intersection $C \cap N$ is the saturation of
$$
I_{C}+I_{N}=(x^2,y,xG,w_0z-z_0w).
$$
Thus, if $G(z_0,w_0)\neq 0$, then $C \cap N$ consists of the reduced point $q$.
Since $G$ has finitely many zeros, corresponding to finitely many points $q_i$ of $L$,
there are finitely lines $N_i=pq_i$ that meets $C$ in a scheme of degree $2$. To finish, it suffices to show that
for each of these lines $N_i$ there is a plane $K$ such that the intersection of $C$ with the divisor
$2N_i$ on $K$ has degree at most three. We claim that, if $q_i=(0,0,z_0,w_0)$, the plane
$K=V(w_0z-z_0w)$ will do. To show this, observe that $w_0z-z_0w$ divides $G$, hence
the ideal of  $C \cap (2N_i)_K$ is the saturation of
$$
(x^2,xy,y^2,w_0z-z_0w, y^{b+1}FP).
$$
This shows $C \cap (2N_i)_K$ is a fat point of degree at most $3$, and completes the proof
that the extremal curves with $g<\frac{1}{2} (d\!-\!2)(d\!-\!3)$
all satisfy (*).

\medskip

This together with Proposition \ref{open-condition} proves half of statemenent \emph{\ref{it:mainThm_a})} of the theorem: if an irreducible component $V$ of $H_{d,g}$ contains an extremal curve then by the openness property, a general curve $C$ in $V$ must also satisfy (*). Thus it remains only to prove \emph{\ref{it:mainThm_b})}, since that statement implies the other half of \emph{\ref{it:mainThm_a})}.

\medskip

\noindent\textbf{Step 3.} We will choose coordinates in
$\PP^3$ adapted to $C$. Let $p$ be a point of $\PP^3$ for which $C$ satisfies (*). Then for every line $N$ containing $p$,
the degree of $ C \cap N$ is at most two, and $ \deg C \cap N = 2$ for at most finitely many lines $N_i$ through $p$.
For each $N_i$, let $K_{ij}$ be the (at most) finitely many planes for which the divisor $2N_{i}$ in $K_{ij}$ meets
$C$ in a scheme of degree $\geq 4$. Finally, choose a line $M$ through $p$ such that a) $M$ does not meet $C$, b)
$M$ is not contained in any plane containing two or more of the $N_i$, and c) $M$ is not contained in any of the planes $K_{ij}$. Choose coordinates $x,y,z,w$ in $\PP^3$ so  that $p=(1,0,0,0)$ and $M$ is the line $z=w=0$.

Let $H$ be the plane $x=0$ and consider the projection $\pi$ of the curve from $p$ to $H$. The image will be a plane curve $C_0 \subset H$ of degree $d$. On the other hand, if we consider the flat family obtained by projecting away from $p$, as in \cite[III, 9.8.3]{AG}, or by applying Proposition \ref{prop:flatWeight} with weight vector $(1,0,0,0)$, the limit is a scheme $C'$ supported on $C_0$, and having nilpotent elements at points where the line $N$ from $p$ meets $C$ twice. Note that the scheme $C'$ is contained in the double plane $2H$ defined by $x^2 = 0$, since no line $N$ meets $C$ more than twice.

In general, when a scheme has embedded points there is no natural scheme structure on the set of embedded points. In our case, however, we can define a scheme $Z$ representing the embedded points by taking the residual of the intersection $C_0$ of $C'$ with $H$ \cite{2h}. Since $C'$ is contained in $2H$, the scheme $Z$ will be contained in $H$ and will have degree equal to the difference of the genus of $C$ (which equals that of $C'$) and of $C_0$, namely $\deg Z = \frac{1}{2}(d-1)(d-2)-g$. One can also obtain the structure sheaf $\mathcal{O}_Z$ as the quotient of $\coo_{C'}$ by $\coo_{C_0}$ or the quotient of $\pi_\ast \coo_{C}$ by $\coo_{C_0}$.

Note that lines $N$ through $p$ meeting $C$ twice correspond to points of $Z$, also in the scheme sense: if $z_0$ is a point of $Z$ of multiplicity $r$, then the cone from $p$ over $z_0$ is an $r$-fold line $N$ meeting $C$ in $2r$ points.
\end{proof}

Before proceeding, we need a lemma.

\bl \label{old-step-one}
Suppose $g \leq \frac{1}{2} (d\!-\!2)(d\!-\!3)$ and let $\nu= \frac{1}{2} (d\!-\!1)(d\!-\!2) - g$.
For any choice of coordinates $x,y,z,w$ on $\PP^3$ and for any curve $C$ in $\PP^3$, the ideal
of $C$ contains a nonzero form $f$ of degree $\nu+1$ of type
\begin{equation}\label{fform}
f(x,y,z,w)= x G_{\nu} (z,w) - \sum_{j=0}^{d-1} y^{j} F_{\nu+1-j}(z,w).
\end{equation}
Here subscripts denote degree.
\el
\begin{proof}
The proof is an easy dimension count.
Inside $H^{0} (\PP^3, \coo_{\PP^3} (\nu+1))$, the subspace $T$ of forms as in \eqref{fform} has dimension
\[
\dim T= \nu+1 + \sum_{j=0}^{d-1} (\nu+2-j)=\nu+1 +d (\nu+2)-\binom{d}{2}= (\nu+1)(d+1)+1 -\frac{1}{2} (d\!-\!1)(d\!-\!2)
\]

We consider the natural map from $H^0\big(\coo_{\PP^3}(\nu+1)\big)$ to $H^0\big(\coo_{C}(\nu+1)\big)$. It is well known -- and easy to show -- that $h^1(\coo_C(n))=0$ if $n \geq d-2$: this is proven in characteristic zero in  \cite[Corollaire 2.4.(1)]{bounds}, and in any characteristic for example in \cite[Proposition 3.2]{sch-beo}.
Now $\nu \geq d-3$ because $g \leq \frac{1}{2} (d\!-\!2)(d\!-\!3)$, hence we can conclude that $h^1(\coo_C(\nu+1))=0$. Therefore by Riemann-Roch
\[
h^0\big( \coo_C (\nu+1)\big) = (\nu+1)d+1-g.
\]
Now we observe that the dimension of $T$ is greater than this number:
\[
\dim T - h^0\big( \coo_C (\nu+1)\big) = \nu+1 - \frac{1}{2}(d-1)(d-2) + g = 1.
\]
Hence there is a non zero form $f\in T$ whose image is zero in $H^0\big(\coo_C(\nu+1)\big)$, that is, which is contained in $H^0\big(\ideal{C}(\nu+1)\big)$.
\end{proof}

\begin{proof}[Proof of Theorem \ref{main-result} (second part)] \textbf{Step 4.} Using the coordinates chosen for $C$ in \textbf{Step 3}, we apply Lemma \ref{old-step-one} and conclude that $C$ is contained in a surface $S$ whose equation has the form \eqref{fform} above. We will be using the surface $S$ for the remainder of the proof. First we will show that $G_\nu$ is not identically zero.

Suppose by way of contradiction $G_\nu = 0$. We can find a plane $K$ with the following properties: it contains the line $M$, it does not contain any of the lines $N_i$ and it is not a component of the surface $S$. We may assume $K$ is the plane of equation $w=0$. Then the curve $S \cap K$ has equation
\[
\sum_{j=0}^{d-1} y^j F_{\nu+1-j}(z,0) = z^{\nu+2-d} \sum_{j=0}^{d-1} c_j\, y^j z^{d-1-j}=0.
\]
Thus, as a divisor in $K$, the intersection $S\cap K$ is a cone consisting of $M$ with multiplicity at least $\nu+2-d$, plus a divisor $D$ of degree at most $d-1$. Since $C$ is contained in $S$, the scheme $C\cap K$ is contained in $S \cap K$, and since $C$ does not meet $M$, this scheme $C \cap K$, which is of degree $d$, must be contained in $D$. Since $K$ contains none of the $N_i$, every line through $p$ in $K$ meets $C \cap K$ at most once, so the divisor $D$ must have degree at least $d$, which is a contradiction.

\noindent\textbf{Step 5.} Now we relate the equation of $S$ to the scheme of embedded points $Z$ from \textbf{Step 3}. Let $q$ be the point $M \cap H$ and let $Z'$ be the projection of $Z$ from $q$ to the line $L$ defined by $x=y=0$. Since $M$ was chosen not to belong to any of the planes $K_{ij}$, it follows that each line through $q$ in $H$ meets $Z$ at most once. Hence the projection from $Z$ to $Z'$ is a closed immersion and $Z'$ is a subscheme of $L$ of degree $\nu$.

On the other hand, the sheaf $\pi_\ast \coo_C$ is generated over $\coo_H$ by the local images of the variable $x$, because we projected away from $p=(1,0,0,0)$. The equation of $S$ shows that
\[
xG_\nu \equiv \sum_{j=0}^{d-1} y^j F_{\nu+1-j}(z,w)\quad \mod I_S.
\]
Thus $G_\nu$ maps the local generator of $\pi_\ast \coo_C$ into $\coo_{C_0}$ and we conclude that $G_\nu$ annihilates $Z$. Since $G_\nu$ and $Z'$ both have the same degree $\nu$, we conclude that $G_\nu (z,w)$ is precisely the annihilator of the scheme $Z'$ in $L$.

\noindent\textbf{Step 6.} Now we will show that $G_\nu$ and $F_{\nu+d-2}$ are forms in $z,w$ with no common zeros. Let $(z_0,w_0)$ be a zero of $G_\nu$. By a linear change of variables we may assume that $(z_0,w_0)=(1,0)$.
 Let $K$ be the plane $w=0$.  Since $K$ contains the line $M$, we conclude from the choices above that $K$ contains at most one of the lines $N_i$, and that the divisor $2N_i$ in $K$ meets $C$ in a scheme of degree $\leq 3$. The plane $K$ is not a component of $S$: otherwise, we would have $S = K+S_1$ where $S_1$ has an equation of the form
 \begin{equation}
 f_1(x,y,z,w) = xG'(z,w) - \sum_{j=0}^{d-1} y^j F'_{\nu-j}(z,w)
 \end{equation}
 where $G'$ is not identically zero, and $\deg G' = \deg G_\nu - 1 = \nu - 1$. As the curve $C$ meets properly every plane through $M$, it must be contained in $S_1$. Then, as in \textbf{Step 5}, we prove that $G'$ kills $Z'$ which is impossible since $Z'$ has degree $\nu$.

 Consider the curve $S \cap K$ which has equation
 \begin{equation}\label{eq_intersec}
\overline{f}=f(x,y,z,0)= - \sum_{j=0}^{d-1} y^{j} F_{\nu+1-j}(z,w)=-z^{\nu+2-d} \sum_{j=0}^{d-1}c_jy^j
z^{d-1-j}\end{equation}
since $G_\nu (z,0)=0$. Now $x$ does not appear in this equation so $S \cap K$ is a divisor in $K$ supported on lines through $p$.
In particular, $\overline{f}$ has a factor $z^{\nu+2-d}$, so the line $M$ appears in $S \cap K$ with multiplicity $\geq \nu+2-d$.
We will show that this multiplicity is exactly $\nu+2-d$, in other words that $c_{d-1} \neq 0$, hence $F_{\nu+2-d}(z,0)$ is not zero, hence
that  $(z,0)$ is not a common zero of $G_\nu$ and  $F_{\nu+2-d}$, as desired.

Since $C$ does not meet the line $M$, the plane $K$ does not contain any component of $C$. Let $W$ be the zero dimensional scheme $C \cap K$. Since $C \subset S$, we find $W \subset S \cap K$. So the curve $S \cap K$ contains
the least divisor $D$ of (possibly multiple) lines through $p$ containing $W$. By construction, most of the reduced lines $N$ through $p$ meet $W$ with degree $1$; at most
one $N_i$ meets $W$ in degree $2$, and for that one, the divisor $2N_i$ in $K$ meets $W$ with degree $\leq 3$. We will show that this divisor $D$ has degree $\geq d-1$, hence $M$ can appear in $S \cap K$ with multiplicity at most
$(\nu+1)-(d-1)=\nu+2-d$ as required.

To verify the assertion that $\deg D \geq d-1$ we argue as follows. If $q\in W$ is a point for which the line $N=pq$
meets $W$ with multiplicity one, then the scheme $W$ is curvilinear at that point, transverse to $N$, so that if its multiplicity at $q$ is $r$, we need $N$ to appear with multiplicity $r$ in $D$ for $D$ to contain $W$.
If $N$ is a line through $p$ meeting $W$ with multiplicity 2, and if $N$ meets $W$ in two distinct points $q_1$,$q_2$,
then, since $2N$ meets $W$ with multiplicity $\leq 3$, one of these points must be a reduced point of $W$, the other is curvilinear transversal to $N$ as above, so that if the multiplicity of the second point is $r$, we need
$N$ to appear with multiplicity $r$ in $D$ for $D$ to contain $W$ at the points $q_1$, $q_2$. Finally, if $N$ meets
$W$ with multiplicity $2$ at a single point $q$, and $2N$ meets it with multiplicity $2$ also, then $\deg W=2$ at $q$, with tangent direction $N$.
The last case is if $N$ meets $W$ at $q$ with multiplicity $2$, and $2N$ meets $W$ at $q$ with multiplicity $3$. This is the only case that requires some calculation, and Lemma \ref{below} below shows in this case that if the multiplicity of
$W$ at $q$ is $r$, then we need $N$ to appear with multiplicity $r\!-\!1$ in $D$ for $D$ to contain $W$ at $q$.
The conclusion is that since $\deg W = d$, we need a divisor $D$ of degree
degree $d-1$ or $d$ to contain $W$, as claimed.
\end{proof}
\renewcommand{\qedsymbol}{$\square$}

\bl \label{below}
Let $A=k[[x,y]]$ with maximal ideal $\mx=(x,y)$. Let $I$ be an $\mx$-primary ideal such that
$\llength \; A/I+(y) = 2$ and $\llength \; A/I+(y^2) = 3$.
Then, replacing if necessary $x$ with another element of
$\mx \setminus \mx^2$, linearly independent of $y$ modulo $\mx^2$, we can put the ideal $I$ in one
of the following forms:
\begin{itemize}
  \item $I=(y-x^2,x^3)$, of colength $3$;
  \item $I=(x^2,xy,y^s)$, with $s \geq 2$, of colength $s+1$
  \item $I=(x^2-y^s,xy,y^{s+1} )$ with $s \geq 2$, of colength $s+2$
\end{itemize}
In all three cases, the least power of $y$ contained in
$I$ is one less than the colength.
\el

\begin{proof}
Suppose $I$ contains an element beginning with a linear form. Since  $A/I+(y)$ has length
$2$ (not $1$), that form must be $y$. Thus $I+(y)=(y,x^2)$, and $I$ contains elements
$y+f$, $x^2+gy$ with $f \in \mx^2$ and $g \in A$. Hence $I$ contains also $x^2-fg \in \mx^2$.
If this element is non-zero, then $I+(y^2)=(y,x^2)$ would have colength only $2$ (not $3$).
Hence $g$ is a unit and $f$ is a unit times $x^2$. We can absorbe the unit into $x$, and obtain the first case
$I=(y-x^2,x^3)$ above.

Suppose now $I \subseteq \mx^2$.   Then
$\length \;A/I+(y^2) = 3 $
implies $\mx^2=I+(y^2)$.
Then $xy-by^2 \in I$ and replacing $x$ with $x-by$ we may assume $xy \in I$.
Since $x^2 \in I + (y^2)$, we can find  $a \in A$
such that  $x^2-ay^2$ belong to $I$. As $xy \in I$,  we can choose such an $a$ in  $k[[y]]$.
If $a=0$, then $I=( x^2,xy,y^{s})$ for some $s \geq 2$. This is the second case above.

If $a \neq 0$, then $ay^2=uy^s$ with $u$ a unit and $s \geq 2$. Then
\[
y^{s+1}=-u^{-1} y(x-uy^s)+u^{-1}xy \in I
\]
In particular we can assume $u$ is a nonzero element in $k$,
and rescaling we can reduce to the case $u=1$.
This gives the third case above (or the second case again if
$y^t \in I$ for some $t \leq s$).
\end{proof}

\begin{proof}[Proof of Theorem \ref{main-result} (third part)]
\textbf{Step 7.} We can now complete the proof of part \emph{\ref{it:mainThm_b})} of the theorem. Let $C$ be a curve satisfying (*). Choose coordinates as in \textbf{Step 3} and a surface $S$ containing $C$ as in \textbf{Step 4} above. We have shown in \textbf{Step 6} that $G_\nu$ and $F_{\nu+d-2}$ have no common zeros. Thus $C$ and $S$ satisfy the hypotheses of Proposition \ref{method}, and we conclude that the flat family constructed there has for its limit an extremal curve. Thus $C$ specializes to an extremal curve, as required.
\end{proof}

\bc \label{smooth}
Every smooth {\em irreducible} curve can be specialized to an extremal curve.
\ec

\begin{proof}
Any smooth irreducible curve satisfies (*). Just choose a point $p \in \PP^3$
that lies on no trisecant line and no tangent line \cite[IV 3.10]{AG}.
\end{proof}

\bc \label{biliaison}
Suppose $C_0$ is a curve satisfying condition (*) that is  minimal in its biliaison class.
Then every curve of degree $d$ and genus $g$
in the biliaison class of $C_0$ is in the connected component of $H_{d,g}$ that contains
the extremal curves of degree $d$ and genus $g$.
\ec

\begin{proof}
The statement follows from Theorem \ref{main-result}.b and \cite[Theorem 2.3]{footnote}.
\end{proof}

\section{Geometric characterization of curves satisfying condition (*)}
The main result of this section is Theorem \ref{geom-condition} that, at least in characteristic zero,
gives a clear geometric characterization of curves that satisfy condition (*).


\bl \label{ed2}
Suppose $C$ satisfies condition (*). Then $C$ has embedding dimension $\leq 2$ at all but finitely many points.
\el

\begin{proof}
If $C$ has embedding dimension $3$ at a point $q$, then every line through $q$ meets $C$ with multiplicity
at least $2$. Since $C$ satisfies (*), this can happen only at finitely many points.
\end{proof}

\bl \label{ed3}
Suppose $C$ is a curve that satisfies condition (*) and has embedding dimension $3$ at a point $q$.
Let $T$ denote the tangent cone to $C$ at $q$: it is naturally embedded in the affine tangent space $\mathbb{A}^3$ of $C$
at $q$. Then one of the following possibilities occurs:
\ben
\item
the cone $T$ is the union of two curves $T_1$ and $T_2$, where $T_1$ is contained in a plane $H \subset \mathbb{A}^3$,
and $T_2$ is a line transversal to the tangent plane to $C_1$;
\item
there is a plane $H$ containing the support of $T$, and the residual curve to $T \cap H$
in $T$ is a reduced line;
\item
the curve $T$ is contained in the complete intersection of two quadric cones in $\mathbb{A}^3$, and therefore has multiplicity
at most $4$ at $q$.
\een

\el

\begin{proof}
Let $p$ be a point with respect to which $C$ satisfies condition (*).
Since $C$ has embedding dimension $3$ at $q$, the line $N=pq$ meets $C$ at $q$ with multiplicity $2$, and
for a general plane $K$ through $N$ the divisor $2N$ in $K$ meets $C$ at $q$ with multiplicity $3$. Fix such a plane
$K$. Then near $q$ the plane section $C \cap K$ has one of the forms listed in Lemma \ref{below}. In the first case
$C \cap K$ has embedding dimension $1$ at $q$, so $C$ has at most embedding dimension $2$ at $q$, a contradiction.
So either the second case or the third case occur. In both cases, the ideal of $C \cap K$ in $\coo_{K,q}$ is contained
in $m^2_{K,q}$, and its image in $m^2/m^3$ has dimension at least $2$. Thus the ideal of $C$ at $q$ also
has image of dimension at least $2$ in  $m_{\PP^3,q}^2/m_{\PP^3,q}^3$, and so the tangent cone is contained
in two distinct quadric cones. The statement then follows immediately.
\end{proof}

\bl \label{nonreduced}
Suppose $C$ satisfies condition (*). Then all non-reduced components of $C$ have support in a single plane.
\el

\begin{proof}
\noindent \textbf{Step 1.}
Suppose $C$ is a non-reduced irreducible curve such that $C_{red}$ is not
a plane curve. We claim that $C$ does not satisfy (*).
Since $C_{red}$ is not
a plane curve, the secant variety of $C_{red}$, that is the closure of the union of lines meeting $C_{red}$ in at least two distinct points,
is the whole of $\PP^3$. Thus a general point $p$ of $\PP^3$ lies on a line
$N=\overline{qr}$ where $q$ and $r$ are distint points of $C_{red}$.
Since $C$ is non-reduced, any plane $H$ through $N$ will intersect $C$ at $q$
in a scheme of length at least two, and therefore the divisor $2N$ on $H$ will intersect $C$ with multiplicity
at least $2$ at $q$. For the same reason,  $2N$  will intersect $C$ with multiplicity
at least $2$ at $r$. Thus $2N$ will meet $C$ at least in a scheme of degree $4$, and
$C$ does not satisfy (*).

\noindent \textbf{Step 2.}
Suppose $C$ and $D$ are  non-reduced irreducible curves that $C_{red} \cup D_{red}$ is not
a plane curve. The $C \cup D$ does not satisfy (*).
Consider the join $J$ of $C_{red}$ and $D_{red}$, that is, the closure of the set of points lying on lines
joining a point of $C_{red}\setminus D_{red}$ with a point of $D_{red}\setminus C_{red}$.
Since $C_{red} \cup D_{red}$ is not
a plane curve, the join $J$ is the whole space $\PP^3$. Thus a general point $p$ of $\PP^3$ lies on a line
$N=\overline{qr}$ with $q \in C_{red}$ and $ r \in D_{red} \setminus C_{red}$.
We may also assume $p$ lies on no plane containing either $C_{red}$ or $D_{red}$.
Then an argument similar to the one above shows that
$C \cup D$ does not satisfy (*).

\noindent \textbf{Step 3.}
Now suppose $C$ is any curve, and that $C_1$, $C_2$, $\ldots C_r$ are the non-reduced components of $C$.
If $C$ satisfies condition (*), then the support of each $C_i$ is contained in a plane $H_i$ by Step 1.
If one of the $C_i$, say $C_1$, is not a line, then $H_1$ is the unique plane containing the support
of $C_1$, and $H_1$ contains the support of the other non-reduced components by Step 2.
Suppose finally that for every $i$ the support of $C_i$ is a line $L_i$. By Step $2$ there is a plane $H$ containing
$L_1$ and $L_2$, and $L_1$ and $L_2$ meet in a point $q$. Suppose by way of contradiction that there is an $i$, say $i=3$, such that $L_i$ is not contained in $H$. By Step 2, the line $L_3$ meets both $L_1$ and $L_2$, and this can happen only at $q$. In particular, at $q$  the curve $C$
has embedding dimension $3$. The tangent cone to $C$ at $q$ then contains at least three (planar) double lines supported on $L_1$, $L_2$ and $L_3$;
but this contradicts Lemma \ref{nonreduced}. Therefore the support of $C$ is contained in the plane $H$ and the statement is proven.
\end{proof}

\bl \label{non-red-comp}
Suppose $C$ is a non-reduced, irreducible curve that satisfies (*).
Then either $C$ has multiplicity $\leq 3$ along its length, or else its underlying multiplicity $2$ structure is planar.
\el

\begin{proof}
Let $D=C_{red}$. The conormal sheaf $\ideal{D}/\ideal{D}^2$ is generically locally free of rank $2$.
Since $C$ has generically embedding dimension $2$,  it does not contain the first infinitesimal neighborhood $D^{(1)}$ of $D$,
that is, $\ideal{C}$ is not contained in $\ideal{D}^2$. Therefore $(\ideal{D}^2+\ideal{C})/\ideal{D}^2$ is generically
a rank one locally free sheaf on $D$, and there is a unique locally Cohen-Macaulay curve $D_2 \subset D^{(1)}$ whose ideal in $D^{(1)}$
coincides with  $(\ideal{D}^2+\ideal{C})/\ideal{D}^2$ at the generic point. This curve $D_2$ is the  multiplicity $2$
structure underlying $C$.  At most points $q$ of $D$, the embedding dimension of $D_2$ is $2$, and there is a well defined
tangent plane $H_q$ to $D_2$. This gives a rational map from $D$ to the dual projective space. If this map is not constant,
then a general point $p$ of $\PP^3$ is contained is some plane $H_q$. If $N=pq$, the  divisor $2N$ in any plane $K$ containing
$N$ meets $C$ with multiplicity $\mu$ at $q$, where $\mu$ is the multiplicity of $C$ along $D$. If $\mu \geq 4$, we thus get
a contradiction to condition (*); thus, if $\mu \geq 4$, $H_q$ is a fixed plane $H$ for every $q$, and $D_2$ is contained in $H$.
\end{proof}

\bt \label{geom-condition}
A curve $C \subset \PP^3$ satisfies condition (*) if and only if it has the following properties:
\begin{enumerate}
  \item[i)] \label{i} there are at most finitely points $q$ at which the embedding dimension of $C$ is $3$;
  at each of these points the embedded tangent cone is contained in a pencil of two dimensional quadratic cones;
  \item[ii)] all non-reduced components of $C$ have support in a single plane;
  \item[iii)] if $C_0$ is a non-reduced component of $C$, then either $C_0$ has multiplicity $\leq 3$ along its length, or else its underlying multiplicity $2$ structure is planar;
  \item[iv)] suppose $C_0$ is a non-reduced component of $C$ that has multiplicity $3$ along its length and whose embedded tangent space at a general point $q \in C_0$ is a plane $H_q$; if $H_q$ varies with $q$,  then we require that  all but finitely many lines $N$ in  $H_q$ meet $C_0$ at $q$ with multiplicity $\leq 2$;
  \item[v)] let $Z \subset \PP^3$ denote the union of the lines that meet $C$ in three or more distinct points; we require the closure of $Z$ not to be the whole of $\PP^3$.
\end{enumerate}
\et

\br
If  $\mbox{char} \, k=0$, conditions {\em iv)} and {\em v)} are automatically satisfied; indeed,
\cite[Corollary 4.6.17]{flenner} implies condition {\em v)}, and {\em iv)} follows from \cite[Theorem 2.1]{hgenus}.
On the other hand, if $k$ has positive characteristic, condition {\em iv)} may fail for curves that satisfy
{\em i)-iii)}. The paper \cite{hgenus} by the first author contains an example. Take any prime $p$, and let $q$ be either $p$ if $p \leq 3$ or 4 if p = 2. Then consider the surface $X$ of equation
$
xw^q + y^3w^{q-2} + yz^q = 0
$.
Let $L$ be the line $x=y=0$, which is contained in the smooth locus of $X$. Take $C = 3L$ on $X$. Then one verifies that the tangent plane of
$C$ varies along the line $L$, and at every point any line in the tangent plane meets the curve $3$ times. As for condition {\em v)},
we asked Joel Roberts about strange curves, and
he gave us an example of a curve in characteristic $p >2$ that is integral, but whose trisecants fill up the whole space.
For $p \geq 3$, take the curve given by affine parametric equation $(t, t^p, t^{p^2})$.
Any secant line meets the curve in $p$ distinct points; thus the $p$-secant lines fill up the whole space. This curve does
not satisfy condition {\em i)} because it has a bad singularity at infinity. We do not know if there is a curve that
does satisfy  {\em i)-iii)} but not {\em v)}.
\er

\begin{proof}
Suppose $C$ satisfies condition (*) with respect to a general point $p$. Then there are at most finitely lines $N$ through $p$ such that $\deg (C \cap N) >1$.
Together with the previous lemmas this property implies  {\em i)-v)}.

If $C$ satisfies {\em i)-v)}, then the set of points $p$  that satisfy the following conditions is open and non-empty:
\begin{enumerate}
  \item \label{1}
   $p$ does not lie on any line that meets $C$ in three or more distinct points;
    \item \label{2}
  $p$ lies on at most finitely many lines meeting $C$ in two distinct points;
    \item \label{3}
  if $q$ is a singular point of $C_{red}$, then the line $pq$ meets $C$ only at $q$
  (there are at most finitely many such points $q$); furthermore, if
  $C$ has embedding dimension two at $q \in C_{red}^{sing}$, then $p$ does not
  belong to the tangent plane to $C$ at $q$;
    \item \label{4}
  $p$ does not lie on any tangent line to $C$ at a smooth point;
    \item \label{5}
    if $q$ is one of the at most finitely many points at which $C$ has embedding dimension three, then
  (a) the line $pq$ meets $C$ only at $q$, (b)  $p$ does not belong
  to the embedded tangent cone $\mathcal{C}_q$ to $C$ at $q$, (c)if the pencil of quadratic cones through $\mathcal{C}_q$
  has a plane $H$ as a common component, then $p \notin H$, and (d) if  $\mathcal{C}_q$ is the complete intersection
  of two quadratic cones, there is a plane $K$ containing $N=pq$ such that $\deg(C \cap (2N)_K) \leq 3$
  (in case (d), if $K$ is a plane meeting $\mathcal{C}_q$ properly at $q$, then for all but finitely lines $N \subset K$ through $q$
  one has $\deg(C \cap (2N)_K) \leq 3$ because $\mathcal{C}_q \cap K$ is a complete intersection of two conics);
  \item \label{6}
  if $C$ has a non-reduced component $C_0$, there are at most finitely many lines through $p$ meeting $C$ at a point
  $q_1 \in C_0$ and at a distinct point $q_2 \in C$; $C$ has embedding dimension $2$ at $q_1$, and the line $pq_1$ does not belong
  to the tangent plane to $C$ at $q_1$; $q_2$ is a smooth point of $C$ (for the support of the non-reduced components of $C$ is contained in a single plane);
    \item \label{7}
  if $q$ is a point at which $C$ has embedding dimension two and the tangent
  plane $H_q$ to $C$ at $q$ contains a non-reduced sub-curve of $C$, then $p$ does not belong $H_q$
  (usually $H_q$ will be the single plane $H$ containing the support of the non-reduced components of $C$;
  let $C_0$ be a non-reduced component of $C$, and let $D \subset H$ be the support of $C_0$; if $D$ is not a line,
  then we must have $H_q=H$; on the other hand, if $D$ is a line, then it is possible that $H_q$ be different from $H$,
  but it will be the only plane containing the planar multiplicity two structure underlying $C_0$, and so it is
  uniquely determined by $C_0$, and there are at most finitely many such planes);
  \item \label{8}
  suppose $C_0$ is a non-reduced component of $C$ that has multiplicity $\leq 3$ along its length;
  let $U$ be as in {\em iv)}; then, if the tangent plane $H_q$ is independent of $q \in U$, $p \notin T_q \; C_0$;
  otherwise, for every $q \in U$ the point $p$ does not belong to any line $N$ through $q$ in  $H_q$ meeting
  $C_0$ at $q$ with multiplicity $\geq 3$.
\end{enumerate}
Choose a point $p$ that satisfies all of the conditions above. We claim that $C$ satisfies condition (*) with respect to $p$.
In the following we denote by $N$ a line through $p$. We have to show: (a) no such line $N$ meets $C$ in a scheme of degree $\geq 3$;
(b) there are at most finitely such lines $N$  that meet $C$ in a scheme of degree $2$, and for each of these lines $N$ there is a plane $K$
containing $N$ such that  $\deg(C \cap (2N)_K) \leq 3$.

By (\ref{1}) any line $N$ through $p$ meets $C$ in at most two distinct points, and by (\ref{2}) at most finitely many of such lines $N$ meet $C$ in two distinct points.
Suppose  $N$ meets $C$ at two distinct points $q_1$ and $q_2$. If $q_1$ and $q_2$ are smooth points of $C$, then $N$ is not tangent to $C$ at these points,
and so $\deg(C \cap N)= 2$; furthermore, if $K$ is a plane through $N$ that does not contain  the tangent line to $C$ neither at $q_1$
nor at $q_2$, then $\deg(C \cap (2N)_K)=2$. Otherwise, we are in the situation of (\ref{6}): $q_1$ is a point of a non-reduced component
$C_0$, $C$ has embedding dimension $2$ at $q_1$, the line $N=pq_1$ does not belong  to the tangent plane to $C$ at $q_1$, and
$q_2$ is a smooth point of $C$. Then $\deg(C \cap N)=2$, and, if $K$ is a plane through $N$ that does not contain  the tangent
line to $C$ at $q_2$, the degree of  $C \cap (2N)_K$ is at most  $3$. Thus we have checked the required properties for lines meeting $C$
in two distinct points.

Suppose now $N$ meets $C$ at a unique point $q$. If $q$ is a smooth point of $C$ or a singular point
of $C$ at which $C$ has embedding dimension $2$ and $N$ is not contained in the tangent plane  $H_q$ to $C$ at $q$, then
$\deg(C \cap N) =1$. Thus we are left to examine lines $N=pq$ where $q$ is either a point such that $C$ has embedding
dimension $2$ at $q$ and  $p$ belongs to the tangent plane $H_q$, or a point at which $C$ has embedding
dimension $3$.  In the first case, by (\ref{3}) $q$ is a smooth point of $C_{red}$, and by (7) it is
contained in a unique non-reduced component $C_0$ whose underlying multiplicity $2$ structure is not planar, so that $C_0$
has multiplicity $\leq 3$ along its length by {\em iii)}; by (\ref{8}) the tangent plane $H_q$ varies with
$q$, so there are finitely many such $q \in C_0$ for which $p \in H_q$; furthermore,  by (\ref{8}) the line $N=pq$ meets
$C$ with multiplicity $2$ at $q$, and for every plane $K$ through $N$, $C \cap (2N)_K$ has degree at most $3$ because
has multiplicity $\leq 3$ along its length.
Finally, suppose $q$ is a point at which $C$ has embedding
dimension $3$. If the tangent cone $\mathcal{C}_q$ is contained in $H \cup L$, where $H$ is a plane and $L$ is a line
transversal to $N$, then $p \notin H$ by (\ref{5}), the line $N=pq$ meets $C$ with multiplicity $2$ at $q$, and, if
$K$ is a plane through $N$ that does not contain $L$, then $(2N)_K$ meets $\mathcal{C}_q$ with multiplicity at most $3$
at $q$, hence $\deg (C \cap (2N)_K) \leq 3$. If the tangent cone $\mathcal{C}_q$ is contained in a plane $H$ plus an embedded
line $L$ , that is, the ideal of $\mathcal{C}_q$ contains $x^2$ and $xy$ up to a choice of coordinates, then
$p \notin H$ by (\ref{5}), the line $N=pq$ meets $C$ with multiplicity $2$ at $q$, and, if
$K$ is a plane through $N$ that does not contain $L$,  $(2N)_K$ meets $C$ with multiplicity at most $3$
at $q$. Finally, the case in which the tangent cone $\mathcal{C}_q$ is the complete intersection of two quadratic cones is taken care
by (\ref{5}.d)

\end{proof}

\bex
Every smooth irreducible curve satisfies condition (*) cf. \ref{smooth}. The same holds
for a smooth non connected curve whose trisecants do not fill the whole space: this latter condition is always satisfied
if the curve lies on a quadric surface,  or if $\mbox{char} \, k=0$. In particular, we obtain a stronger version of the main
result of \cite{ls}: every smooth divisor on a smooth quadric surface specializes in a flat family to an extremal curve.
\eex

\bex
Assume $\mbox{char} \, k=0$. Let $C$ be a reduced curve that has embedding dimension at most
$2$ at every point. Then $C$ satisfies condition (*). This applies for example to  general ACM curves:
see \cite[Theorem 7.21]{hs}.
\eex


\bex
Assume $\mbox{char} \, k=0$. Suppose we are given curves $C_i$ for $i=0,\ldots,r$ that are disjoint,
satisfy condition (*) and are reduced except possibly for $C_0$. Then their disjoint union $C$ is a curve
that satisfies condition (*).
\eex

\bex
Assume $\mbox{char} \, k=0$. Suppose  $C$ satisfies condition (*) and $D$ is a smooth irreducible curve
meeting $C$ at a single point $p$ in such a way that the tangent line to $D$ at $p$ is not contained
in the tangent space to $C$ at $p$. Then $C \cup D$ satisfies condition (*).
\eex

\end{document}